\documentclass{amsart}
\usepackage{mathrsfs}
\usepackage{amsfonts}

\newtheorem{theorem}{Theorem}[section]
\newtheorem{lemma}[theorem]{Lemma}

\theoremstyle{definition}

 \theoremstyle{remark}
\newtheorem{remark}[theorem]{Remark}
\newtheorem{corollary}[theorem]{Corollary}

 \numberwithin{equation}{section}

\begin{document}

\title{On the best constants of Hardy inequality in $\mathbb{R}^{n-k}\times (\mathbb{R}_{+})^{k}$ and related improvements}

%    Information for first author
\author{ Dan Su}
%    Address of record for the research reported here
\address{School of Information Technology and Management Engineering,
      \ University of International Business and Economics, Beijing, 100029, People's Republic of China}
\email{sudan@uibe.edu.cn}
\author{Qiao-Hua Yang}
%    Address of record for the research reported here
\address{School of Mathematics and Statistics, Wuhan University, Wuhan, 430072, People's Republic of China}
%    Current address
\curraddr{L.M.A.M., Universit\'e de Bretagne-Sud, Centre de
Recherche, Campus de Tohannic, BP 573, F-56017 Vannes, France }
\email{qhyang.math@gmail.com}
%    \thanks will become a 1st page footnote.
\thanks{This work is supported by Program for  Innovative Research Team in UIBE
and   the National Natural Science Foundation of
China(No.11101096).}

%    Information for second author
%\author{Author Two}
%\address{Mathematical Research Section, School of Mathematical Sciences,
%Australian National University, Canberra ACT 2601, Australia}
%\email{two@maths.univ.edu.au}
%\thanks{Support information for the second author.}

%    General info

%\date{January 1, 2001 and, in revised form, June 22, 2001.}

%\dedicatory{This paper is dedicated to our advisors.}


\begin{abstract}
We compute the explicit sharp constants of  Hardy  inequalities in
the cone $\mathbb{R}_{k_+}^{n}:=\mathbb{R}^{n-k}\times
(\mathbb{R}_{+})^{k}=\{(x_{1},\cdots,x_{n})|x_{n-k+1}>0,\cdots,x_{n}>0\}$
with $1\leq k\leq n$. Furthermore, the spherical harmonic
decomposition is given for a function $u\in
C^{\infty}_{0}(\mathbb{R}_{k_+}^{n})$. Using this decomposition and
following the idea of Tertikas and Zographopoulos, we obtain the
Filippas-Tertikas improvement of the Hardy inequality.
\end{abstract}

\maketitle

%\section*{This is an unnumbered first-level section head}

\section{Introduction}
Let $\Sigma$ be a domain in $\mathbb{S}^{n-1}$, the unit sphere in
$\mathbb{R}^{n}$,  and let $\mathcal{C}_{\Sigma}\subset
\mathbb{R}^{n}$ be the cone associated with $\Sigma$:
\[
\mathcal{C}_{\Sigma}:=\{t\sigma| t>0,\;\;\sigma\in \Sigma\}.
\]
The Hardy inequality in $\mathcal{C}_{\Sigma}$ states that, for all
$u\in C^{\infty}_{0}(\mathcal{C}_{\Sigma})$, there holds (cf.
\cite{n,mm})
\begin{equation}\label{1.1}
\int_{\mathcal{C}_{\Sigma}}|\nabla u(x)|^{2}dx \geq
\left(\frac{(n-2)^{2}}{4}+\lambda_{1}(\Sigma)\right)\int_{\mathcal{C}_{\Sigma}}\frac{u(x)^{2}}{|x|^{2}}dx
\end{equation}
and the constant
$\left(\frac{(n-2)^{2}}{4}+\lambda_{1}(\Sigma)\right)$ in
(\ref{1.1}) is sharp, where $\lambda_{1}(\Sigma)$ is the Dirichlet
principal eigenvalue of the spherical Laplacian
$-\Delta_{\mathbb{S}^{n-1}}$ on $\Sigma$. In some special cases, the
exact value of $\lambda_{1}(\Sigma)$ can be computed. We note the
 value of $\lambda_{1}(\Sigma)$ has been full-filled in the
case of $n=2$ (cf. \cite{cm1}). To the best of our knowledge (cf.
\cite{cm2,c1,c2,m,mm,n}), when $n\geq 3$ , $\lambda_{1}(\Sigma)$ is
known only in the case of $\Sigma=\mathbb{S}_{+}^{n-1}$, the
semi-sphere mapped in the upper half space
$\mathbb{R}^{n}_{+}=\{(x_{1},\cdots,x_{n})|x_{n}>0\}$. In fact, it
can be computed via
  the following sharp Hardy inequality (cf.
\cite{ftt})
\begin{equation}\label{1.2}
\int_{\mathbb{R}^{n}_{+}}|\nabla u(x)|^{2}dx \geq
\frac{n^{2}}{4}\int_{\mathbb{R}^{n}_{+}}\frac{u(x)^{2}}{|x|^{2}}dx.
\end{equation}

One of the aim of this note is to  compute the explicit sharp
constants of  Hardy  inequalities in the cone
$\mathbb{R}_{k_+}^{n}=\{(x_{1},\cdots,x_{n})|x_{n-k+1}>0,\cdots,x_{n}>0\}$,
where $1\leq k\leq n$.   To this end, we have:
\begin{theorem}
Let $n\geq 3$.
 There holds, for all $u\in C^{\infty}_{0}(\mathbb{R}_{k_+}^{n})$,
\begin{equation}\label{1.3}
\int_{\mathbb{R}_{k_+}^{n}}|\nabla u|^{2}dx \geq
\frac{(n-2+2k)^{2}}{4}
\int_{\mathbb{R}_{k_+}^{n}}\frac{u^{2}}{|x|^{2}}dx,
\end{equation}
and the constant $\frac{(n-2+2k)^{2}}{4}$ in (\ref{1.3}) is sharp.
\end{theorem}
We note the proof of Theorem 1.1 above is similar to that of Theorem
1.2 and Corollary 1.3 in \cite{ly} and also to that of Theorem 6.1
in \cite{ms}. Combing the inequality (\ref{1.1}) and Theorem 1.1
yields
\begin{corollary}
$\lambda_{1}(\mathbb{S}^{n-1}\cap \mathbb{R}_{k_+}^{n})=k(n+k-2)$
 for all $n\geq 3$.
\end{corollary}

Next, we consider the spherical harmonic decomposition of a function
$u\in C^{\infty}_{0}(\mathbb{R}_{k_+}^{n})$. We show that for a
function $u\in C^{\infty}_{0}(\mathbb{R}_{k_+}^{n})$, it has the
expansion in spherical harmonics ( for details, see section 3)
\[
u(x)=\sum^{\infty}_{l=k}f_{l}(r)\phi_{l}(\sigma),
\]
where $r=|x|$ and $\phi_{l}(\sigma)$ ($l\geq k$) are the orthonormal
eigenfunctions of the spherical Laplacian
$-\Delta_{\mathbb{S}^{n-1}}$ with responding eigenvalues $l(n+l-2)$.
Using this decomposition and following the idea of Tertikas and
Zographopoulos (\cite{tz}), one can easily obtain several
improvements of inequality (\ref{1.3}) when $u$ is supported in a
bounded domain $\Omega\subset \mathbb{R}_{k_+}^{n}$. For example,
      we have the following
Filippas-Tertikas improvement (cf. \cite{ft}):
\begin{theorem}
Let $n\geq 3$.
 There holds, for all $u\in C^{\infty}_{0}(B_{R}\cap\mathbb{R}_{k_+}^{n})$,
\begin{equation*}
\int_{B_{R}\cap\mathbb{R}_{k_+}^{n}}|\nabla u|^{2} \geq
\frac{(n-2+2k)^{2}}{4}
\int_{B_{R}\cap\mathbb{R}_{k_+}^{n}}\frac{u^{2}}{|x|^{2}}+\frac{1}{4}\sum^{\infty}_{i=1}\int_{B_{R}
\cap\mathbb{R}_{k_+}^{n}}\frac{u^{2}}{|x|^{2}}X^{2}_{1}\left(\frac{|x|}{R}\right)\cdot\cdots
\cdot X^{2}_{i}\left(\frac{|x|}{R}\right),
\end{equation*}
where \[ X_{1}(s)=(1-\ln s)^{-1},\;\;X_{i}(s)=X_{1}(X_{i-1}(t))
\]
for $i\geq 2$ and $B_{R}=\{x\in \mathbb{R}^{n}:|x|<R\}$.
\end{theorem}

\section{Proof of Theorem 1.1}
Let $l>0$. A simple calculation shows, for $x_{n}>0$,
\begin{equation}\label{b1}
(x_{n})^{-l}\left(-\Delta
  +\frac{l(l-1)}{x^{2}_{n}}\right)(x^{l}_{n}g(x))
  =-\left(\sum^{n}_{j=1}\frac{\partial^{2}}{\partial x^{2}_{j}}+\frac{2l}{x_{n}}\frac{\partial}{\partial
  x_{n}}\right)g(x).
\end{equation}
Notice that $\frac{\partial^{2}}{\partial
  x^{2}_{n}}+\frac{2l}{x_{n}}\frac{\partial}{\partial
  x_{n}}$ is nothing but the $(2l+1)$-dimensional Laplacian of a radial
  function if  $2l$ is a positive integer. So following the
proof of Theorem 1.2 in \cite{ly} or Theorem 6.1 in \cite{ms}, we
have:
\begin{lemma}
There holds, for $l\in \{1/2,1,3/2,2,\cdots,n/2,\cdots\}$ and $u\in
C^{\infty}_{0}(\mathbb{R}_{+}^{n})$,
\begin{equation}\label{2.2}
\int_{\mathbb{R}^{n}_{+}}|\nabla u|^{2}dx
+l(l-1)\int_{\mathbb{R}^{n}_{+}}\frac{u^{2}}{x_{n}^{2}}dx\geq
\frac{(n+2l-2)^{2}}{4}\int_{\mathbb{R}^{n}_{+}}\frac{u^{2}}{|x|^{2}}dx
\end{equation}
and the constant $\frac{(n+2l-2)^{2}}{4}$ in (\ref{2.2}) is sharp.
\end{lemma}
\begin{proof}
Recall  the sharp Hardy inequality  on
$\mathbb{R}_{x}^{n-1}\times\mathbb{R}_{y}^{2l+1}$:
\begin{equation}\label{a3.1}
\int_{\mathbb{R}_{x}^{n-1}\times\mathbb{R}_{y}^{2l+1}}|\nabla
v|^{2}\geq\frac{(n+2l-2)^{2}}{4}\int_{\mathbb{R}_{x}^{n-1}\times\mathbb{R}_{y}^{2l+1}}
\frac{v^{2}}{x^{2}_{1}+\cdots+x_{n-1}^{2}+|y|^{2}},
\end{equation}
where  $v\in
C^{\infty}_{0}(\mathbb{R}_{x}^{n-1}\times\mathbb{R}_{y}^{2l+1})$.
The constant that appear in (\ref{a3.1}) is also sharp if one
consider only the functions like $\widetilde{v}(x,|y|)\in
C^{\infty}_{0}(\mathbb{R}_{x}^{n-1}\times\mathbb{R}_{y}^{2l+1})$.
Set $x_{n}=|y|$ and
$\varphi(x_{1},\cdots,x_{n})=\widetilde{v}(x,|y|)$, we can deduce,
by (\ref{a3.1}) and (\ref{b1}),
\begin{equation*}
\begin{split}
\int_{\mathbb{R}_{x}^{n-1}\times\mathbb{R}_{y}^{2l+1}}|\nabla
\widetilde{v}|^{2}&=-\int_{\mathbb{R}_{x}^{n-1}\times\mathbb{R}_{y}^{2l+1}}\widetilde{v}(x,|y|)\left(\sum^{n-1}_{j=1}\frac{\partial^{2}}{\partial
x^{2}_{j}}+\sum^{2l+1}_{k=1}\frac{\partial^{2}}{\partial
y^{2}_{k}}\right)\widetilde{v}(x,|y|)\\
&=-\int_{\mathbb{R}_{x}^{n-1}\times\mathbb{R}_{y}^{2l+1}}\varphi(x)\left(\sum^{n}_{j=1}\frac{\partial^{2}}{\partial
x^{2}_{j}}+\frac{2l}{x_{n}}\frac{\partial}{\partial
  x_{n}}\right)\varphi(x)\nonumber\\
  &=-\int_{\mathbb{R}_{x}^{n-1}\times\mathbb{R}_{y}^{2l+1}}x^{-l}\varphi(x)\left(\sum^{n}_{j=1}\frac{\partial^{2}}{\partial
x^{2}_{j}}+\frac{l(l-1)}{x^{2}_{n}}\right)(x^{l}\varphi(x))\nonumber\\
&=-|\mathbb{S}^{2l+1}|\int_{\mathbb{R}^{n}_{+}}x^{l}\varphi(x)\left(\sum^{n}_{j=1}\frac{\partial^{2}}{\partial
x^{2}_{j}}+\frac{l(l-1)}{x^{2}_{n}}\right)(x^{l}\varphi(x))\\
&\geq\frac{(n+2l-2)^{2}}{4}\int_{\mathbb{R}_{x}^{n-1}\times\mathbb{R}_{y}^{2l+1}}
\frac{\widetilde{v}^{2}}{x^{2}_{1}+\cdots+x_{n-1}^{2}+|y|^{2}}\\
&=\frac{(n+2l-2)^{2}|\mathbb{S}^{2l+1}|}{4}\int_{\mathbb{R}_{+}^{n}}
\frac{\varphi^{2}x^{2l}_{n}}{|x|^{2}},
\end{split}
\end{equation*}
where $|\mathbb{S}^{2l+1}| $ is the volume of $\mathbb{S}^{2l+1}$.
It remains  to set $u=x^{l}_{n}\varphi$.
\end{proof}

\begin{remark}
If we let $l(l-1)=0$ in Lemma 2.1, then  $l=1$ and we obtain the
sharp Hardy inequality on the half space $\mathbb{R}^{n}_{+}$ (see
\cite{ftt} for a different proof)
\[
\int_{\mathbb{R}^{n}_{+}}|\nabla u(x)|^{2}dx \geq
\frac{n^{2}}{4}\int_{\mathbb{R}^{n}_{+}}\frac{u(x)^{2}}{|x|^{2}}dx.
\]
Notice that this inequality is one of the objects in the Theorem 1.1
and  the  dimension $2l+1=3$ play an important role.  So, in order
to prove Theorem 1.1,  we can  repeat the same argument of Corollary
1.3 in \cite{ly} by choosing such  dimension $3$.
\end{remark}

\emph{Proof of Theorem 1.1.} Notice that
\begin{equation}\label{b2}
\begin{split}
&-\prod^{n}_{i=n-k+1}x_{i}^{-1}\sum^{n}_{j=1}\frac{\partial^{2}}{\partial
x^{2}_{j}}
  \left(\prod^{n}_{i=n-k+1}x_{i}g(x)\right)\\
  =&-\sum^{n-k}_{j=1}\frac{\partial^{2}g(x)}{\partial x^{2}_{j}}-
   \sum^{n}_{j=n-k+1}\left(\frac{\partial^{2}}{\partial x^{2}_{j}}+\frac{2}{x_{j}}\frac{\partial}{\partial
  x_{j}}\right)g(x).
\end{split}
\end{equation}
We consider  the sharp Hardy inequality on
$\mathbb{R}_{x}^{n-k}\times\mathbb{R}_{y}^{3k}$:
\begin{equation*}
\int_{\mathbb{R}_{x}^{n-k}\times\mathbb{R}_{y}^{3k}}(|\nabla_{x}
v|^{2}+|\nabla_{y}
v|^{2})\geq\frac{(n+2k-2)^{2}}{4}\int_{\mathbb{R}_{x}^{n-k}\times\mathbb{R}_{y}^{3k}}
\frac{v^{2}}{\sum^{n-k}_{i=1}x^{2}_{i}+\sum^{3k}_{j=1}y^{2}_{j}},
\end{equation*}
where $ v\in
C^{\infty}_{0}(\mathbb{R}_{x}^{n-k}\times\mathbb{R}_{y}^{3k}).$ Set
$$x_{n-k+1}=\sqrt{y^{2}_{1}+y^{2}_{2}+y^{2}_{3}},\;\;x_{n-k+2}=\sqrt{y^{2}_{4}+y^{2}_{5}+y^{2}_{6}},\cdots,
x_{n}=\sqrt{y^{2}_{3k-2}+y^{2}_{3k-1}+y^{2}_{3k}}$$  and consider
all the functions like
$$v(x_{1},\cdots,x_{n-k},y_{1},\cdots,y_{3k})=\widetilde{v}(x_{1},\cdots,x_{n}).$$
The constant $\frac{(n+2k-2)^{2}}{4}$  is also sharp for such
functions (see e.g. \cite{ssw}). Following the proof of Lemma 2.1,
we have, using (\ref{b2}),
\begin{equation*}
\begin{split}
\int_{\mathbb{R}_{x}^{n-1}\times\mathbb{R}_{y}^{3k}}|\nabla
\widetilde{v}|^{2}
&=-|\mathbb{S}^{3}|^{k}\int_{\mathbb{R}^{n}_{k_{+}}}\prod^{n}_{i=n-k+1}x_{i}\widetilde{v}(x)\sum^{n}_{j=1}\frac{\partial^{2}}{\partial
x^{2}_{j}}\left(\prod^{n}_{i=n-k+1}x_{i}\widetilde{v}(x)\right)\\
&\geq\frac{(n+2k-2)^{2}}{4}\int_{\mathbb{R}_{x}^{n-k}\times\mathbb{R}_{y}^{3k}}
\frac{\widetilde{v}^{2}}{\sum^{n-k}_{i=1}x^{2}_{i}+\sum^{3k}_{j=1}y^{2}_{j}}\\
&=\frac{(n+2k-2)^{2}|\mathbb{S}^{3}|^{k}}{4}\int_{\mathbb{R}_{k_{+}}^{n}}
\frac{\widetilde{v}^{2}\prod^{n}_{i=n-k+1}x^{2}_{i}}{|x|^{2}}.
\end{split}
\end{equation*}
It remains  to set $u=\widetilde{v}\prod^{n}_{i=n-k+1}x_{i}$ and the
desired result follows.

\section{spherical harmonic decomposition}
For a function $u\in C^{\infty}_{0}(\mathbb{R}_{k_+}^{n})$, we
denote by $\widetilde{u}$ the odd extension of variables
$\{x_{n-k+1},\cdots,x_{n}\}$ of $u$, i.e. $\widetilde{u}(x)$
satisfies
\[
\widetilde{u}(x_{1},\cdots,x_{n})=u(x_{1},\cdots,x_{n}),\;\;\forall
(x_{1},\cdots,x_{n})\in \mathbb{R}_{k_+}^{n}
\]
and
\[
\widetilde{u}(x_{1},\cdots,x_{j-1},-x_{j},x_{j+1},\cdots,x_{n})=-\widetilde{u}(x_{1},\cdots,x_{j-1},x_{j},x_{j+1},\cdots,x_{n})
\]
for all $n-k+1\leq j\leq n$. Then $\widetilde{u}\in
C^{\infty}_{0}(\mathbb{R}^{n})$ and moreover,
\begin{equation}\label{3.1}
\int_{\mathbb{R}^{n}}|\nabla
\widetilde{u}|^{2}=2^{k}\int_{\mathbb{R}_{k_{+}}^{n}}|\nabla
u|^{2},\;\;\;\;\int_{\mathbb{R}^{n}}\frac{
\widetilde{u}^{2}}{|x|^{2}}=2^{k}\int_{\mathbb{R}_{k_{+}}^{n}}\frac{
u^{2}}{|x|^{2}}.
\end{equation}
Decomposing $\widetilde{u}$ into spherical harmonics we get (see
e.g. \cite{tz})
\begin{equation}\label{3.2}
\widetilde{u}=\sum^{\infty}_{l=0}\widetilde{u}_{l}:=\sum^{\infty}_{l=0}f_{l}(r)\phi_{l}(\sigma),
\end{equation}
 where $\phi_{l}(\sigma)$ are the orthonormal
eigenfunctions of the Laplace-Beltrami operator with responding
eigenvalues
$$c_{l}=l(n+l-2),\;\;l\geq 0.$$
The functions $f_{l}(r)$ belong to $C^{\infty}_{0}(\mathbb{R}^{n})$,
satisfying $f_{l}(r)=O(r^{l})$ and $f'_{l}(r)=O(r^{l-1})$ as
$r\rightarrow0$. Without loss of generality, we  assume
\[\int_{\mathbb{S}^{n-1}}|\phi_{l}(\sigma)|^{2}d\sigma=1,\;\;\forall l\geq 0.\]
By (\ref{3.2}),
\[
f_{l}(r)=\int_{\mathbb{S}^{n-1}}\widetilde{u}(x)\phi_{l}(\sigma)
d\sigma
\]
and
\begin{equation}\label{3.3}
\begin{split}
\int^{\infty}_{0}f^{2}_{l}(r)r^{n+l-1}dr=&\int^{\infty}_{0}\int_{\mathbb{S}^{n-1}}\widetilde{u}(x)f_{l}(r)\phi_{l}(\sigma)
r^{n+l-1}d\sigma
dr\\
=&\int_{\mathbb{R}^{n}}\widetilde{u}(x)f_{l}(|x|)\phi_{l}(\sigma)|x|^{l}dx.
\end{split}
\end{equation}

\begin{lemma}
$f_{l}=0$ for all $0\leq l\leq k-1$.
\end{lemma}
Before the proof of Lemma 3.1, we need some multi-index notation. We
denote by $\mathbb{N}_{0}$ the set of nonnegative integer. A
multi-index is denoted by
$\alpha=(\alpha_{1},\cdots,\alpha_{n})\in\mathbb{N}_{0}^{n}$. For
$\alpha \in \mathbb{N}_{0}^{n}$ and $x\in \mathbb{R}^{n}$ a monomial
in variables $x_{1},\cdots, x_{n}$ of index $\alpha$ is defined by
\[
x^{\alpha}=x^{\alpha_{1}}_{1}\cdots x^{\alpha_{n}}_{n}.\]
 The number
$|\alpha|=\alpha_{1}+\cdots+\alpha_{n}$ is called the total degree
of $x^{\alpha}$. Notice that $\phi_{l}(\sigma)$ is nothing but the
spherical harmonic of degree $l$ (see e.g. \cite{sw}, Chapter IV),
it has the expansion
\begin{equation}\label{3.4}
\phi_{l}(\sigma)=\frac{1}{|x|^{l}}\sum_{|\alpha|=l}C_{\alpha}x^{\alpha}
\end{equation}
for some constants $C_{\alpha}\in\mathbb{R}$.

\emph{Proof of Lemma 3.1.} By (\ref{3.3}) and (\ref{3.4}),
\begin{equation*}
\begin{split}
\int^{\infty}_{0}f^{2}_{l}(r)r^{n+l-1}dr&=\int_{\mathbb{R}^{n}}\widetilde{u}(x)f_{l}(|x|)\phi_{l}(\sigma)|x|^{l}dx\\
&=\sum_{|\alpha|=l}C_{\alpha}\int_{\mathbb{R}^{n}}\widetilde{u}(x)f_{l}(|x|)x^{\alpha}dx.
\end{split}
\end{equation*}
So to finish the proof, it is enough to show
\[
\int_{\mathbb{R}^{n}}\widetilde{u}(x)f_{l}(|x|)x^{\alpha}dx=0
\]
for all $|\alpha|=l$ with $0\leq l\leq k-1$.

For $|\alpha|=\alpha_{1}+\cdots+\alpha_{n}=l\leq k-1$, there must
exist  $j$, $n-k+1\leq j\leq n$, such that $\alpha_{j}=0$ (we note
if $\alpha_{j}>0$ for all $n-k+1\leq j\leq n$, then
$\alpha_{n-k+1}+\cdots+\alpha_{n}\geq k$ and this is a contradiction
to $|\alpha|\leq k-1$). Therefore,
\begin{equation*}\label{3.6}
\begin{split}
&\int_{\mathbb{R}^{n}}\widetilde{u}(x)f_{l}(|x|)x^{\alpha}dx=
\int_{\mathbb{R}^{n}}\widetilde{u}(x)f_{l}(|x|)
x^{\alpha_{1}}_{1}\cdots x^{\alpha_{j-1}}_{j-1}\cdot
x^{\alpha_{j+1}}_{j+1}\cdots x^{\alpha_{n}}_{n}dx\\
=&\int_{\mathbb{R}^{n-1}}x^{\alpha_{1}}_{1}\cdots
x^{\alpha_{j-1}}_{j-1}\left(\int_{\mathbb{R}}\widetilde{u}(x)f_{l}(|x|)dx_{j}\right)
x^{\alpha_{j+1}}_{j+1}\cdots x^{\alpha_{n}}_{n}dx_{1}\cdots
dx_{j-1}dx_{j+1}\cdots dx_{n}.
\end{split}
\end{equation*}
Since $\widetilde{u}(x)$ is an odd function of variable $x_{j}$, so
does $\widetilde{u}(x)f_{l}(|x|)$. Therefore,
\[\int_{\mathbb{R}}\widetilde{u}(x)f_{l}(|x|)dx_{j}=0\]
and hence
\[\int_{\mathbb{R}^{n}}\widetilde{u}(x)f_{l}(|x|)x^{\alpha}dx=0.\]
The proof of Lemma 3.1 is now completed.

\begin{remark}
By Lemma 3.1, the function $\widetilde{u}$,  \;the odd extension of
variables $\{x_{n-k+1},\cdots,x_{n}\}$ of $u$,  has the expansion in
spherical harmonics
\[
u(x)=\sum^{\infty}_{l=k}f_{l}(r)\phi_{l}(\sigma),
\]
so does the function $u$ itself in $\mathbb{R}_{k_+}^{n}$.
\end{remark}

\emph{Proof of Theorem 1.3.}  If we extend $u$ as zero in
$\mathbb{R}_{k_+}^{n}\setminus B_{R}$, we may consider $u\in
C^{\infty}_{0}(\mathbb{R}_{k_+}^{n})$. By (\ref{3.1}), it is enough
to show
\[
\int_{B_{R}}|\nabla \widetilde{u}|^{2} \geq \frac{(n-2+2k)^{2}}{4}
\int_{B_{R}}\frac{\widetilde{u}^{2}}{|x|^{2}}+\frac{1}{4}\sum^{\infty}_{i=1}\int_{B_{R}
}\frac{\widetilde{u}^{2}}{|x|^{2}}X^{2}_{1}\left(\frac{|x|}{R}\right)\cdot\cdots
\cdot X^{2}_{i}\left(\frac{|x|}{R}\right)
\]
hold for all $\widetilde{u}\in C^{\infty}_{0}(B_{R})$. Since
$\widetilde{u}$ has the expansion in spherical harmonics
\[
u(x)=\sum^{\infty}_{l=k}f_{l}(r)\phi_{l}(\sigma),
\]
where $f_{l}(r)\in C^{\infty}_{0}(B_{R})$, satisfying
$f_{l}(r)=O(r^{l})$ and $f'_{l}(r)=O(r^{l-1})$ as $r\rightarrow0$,
we have,
\begin{equation*}
\begin{split}
&\int_{B_{R}}|\nabla \widetilde{u}|^{2} -\frac{(n-2+2k)^{2}}{4}
\int_{B_{R}}\frac{\widetilde{u}^{2}}{|x|^{2}}\\
=&\sum^{\infty}_{l=k}\left[\int_{B_{R}}|f'_{l}(r)|^{2}dx
+l(n+l-2)\int_{B_{R}}\frac{f_{l}^{2}(r)}{|x|^{2}}dx-\frac{(n-2+2k)^{2}}{4}
\int_{B_{R}}\frac{f_{l}^{2}(r)}{|x|^{2}}dx\right]\\
=&\sum^{\infty}_{l=k}\left[\int_{B_{R}}|f'_{l}(r)|^{2}dx
+(l-k)(n+l+k-2)\int_{B_{R}}\frac{f_{l}^{2}(r)}{|x|^{2}}dx-\frac{(n-2)^{2}}{4}
\int_{B_{R}}\frac{f_{l}^{2}(r)}{|x|^{2}}dx\right]\\
\geq&\sum^{\infty}_{l=k}\left[\int_{B_{R}}|f'_{l}(r)|^{2}dx
-\frac{(n-2)^{2}}{4}
\int_{B_{R}}\frac{f_{l}^{2}(r)}{|x|^{2}}dx\right].
\end{split}
\end{equation*}
To get the last inequality above, we use the fact
$(l-k)(n+l+k-2)\geq0$ since $l\geq k\geq 1$. Recalling the
Filippas-Tertikas improvement of Hardy inequality (cf. \cite{ft,tz})
\[
\int_{B_{R}}|f'_{l}(r)|^{2}dx -\frac{(n-2)^{2}}{4}
\int_{B_{R}}\frac{f_{l}^{2}(r)}{|x|^{2}}dx\geq
\frac{1}{4}\sum^{\infty}_{i=1}\int_{B_{R}
}\frac{f_{l}^{2}(r)}{|x|^{2}}X^{2}_{1}\left(\frac{|x|}{R}\right)\cdot\cdots
\cdot X^{2}_{i}\left(\frac{|x|}{R}\right),
\]
we have
\begin{equation*}
\begin{split}
\int_{B_{R}}|\nabla \widetilde{u}|^{2} -\frac{(n-2+2k)^{2}}{4}
\int_{B_{R}}\frac{\widetilde{u}^{2}}{|x|^{2}}
\geq&\sum^{\infty}_{l=k}\left[\int_{B_{R}}|f'_{l}(r)|^{2}dx
-\frac{(n-2)^{2}}{4}
\int_{B_{R}}\frac{f_{l}^{2}(r)}{|x|^{2}}dx\right]
\\
\geq&\frac{1}{4}\sum^{\infty}_{l=k}\int_{B_{R}
}\frac{f_{l}^{2}(r)}{|x|^{2}}X^{2}_{1}\left(\frac{|x|}{R}\right)\cdot\cdots
\cdot X^{2}_{i}\left(\frac{|x|}{R}\right)\\
=&\frac{1}{4}\int_{B_{R}
}\frac{\widetilde{u}^{2}(r)}{|x|^{2}}X^{2}_{1}\left(\frac{|x|}{R}\right)\cdot\cdots
\cdot X^{2}_{i}\left(\frac{|x|}{R}\right).
\end{split}
\end{equation*}
The desired result follows.
\section*{Acknowledgements}
The authors thanks the referee for his/her careful reading and very
useful comments which improved the final version of this paper.

%\section*{Acknowledgements}
%The authors thanks the referee for his/her careful reading and very
%useful comments which improved the final version of this paper.

\end{document}